\documentclass{amsart}
\usepackage{amsmath} 
\usepackage{amssymb}
\usepackage{mathrsfs}

\newtheorem{theorem}{Theorem}[section]

\newtheorem{proposition}[theorem]{Proposition}

\theoremstyle{definition}
\newtheorem{definition}[theorem]{Definition}

\newtheorem{discussion}[theorem]{Discussion}

\theoremstyle{remark}
\newtheorem{remark}[theorem]{Remark}

\numberwithin{equation}{section}
\numberwithin{figure}{section}
\newtheorem{thm}{Theorem}[section]
\newtheorem{defn}[thm]{Definition}

\newtheorem{conv}[thm]{Convention}
\newtheorem{disc}[thm]{Discussion}

\newcommand{\lf}{{\rm lf}}

\newcommand{\tp}{{\rm tp}}

\newcommand{\cl}{{\rm cl}}

\newcommand{\Ex}{{\rm Ex}}

\newcommand{\eq}{{\rm eq}}

\newcommand{\EC}{{\rm EC}}

\newcommand{\Dom}{{\rm Dom}}

\newcommand{\wilog}{{\rm without loss of generality}}
\newcommand{\Wilog}{{\rm Without loss of generality}}

\newcommand{\mn}{{\medskip\noindent}}
\newcommand{\sn}{{\smallskip\noindent}}

\newcommand{\gC}{{\mathfrak C}}

\newcommand{\cI}{{\mathscr I}}

\newcommand{\gk}{{\mathfrak k}}

\newcommand{\bS}{{\bf S}}

\newcommand{\gt}{{\mathfrak t}}

\newcommand{\M}{{\mathcal{M}}}	

\newcommand{\I}{{\bf I}}	
\newcommand{\J}{{\bf J}}

\renewcommand{\lg}{{\ell\!g}}

\newcommand{\mfk}{{\mathfrak{k}}}
\newcommand{\keq}{{\mfk^{\rm eq}}}

\newcount\skewfactor
\def\mathunderaccent#1#2 {\let\theaccent#1\skewfactor#2
\mathpalette\putaccentunder}
\def\putaccentunder#1#2{\oalign{$#1#2$\crcr\hidewidth
\vbox to.2ex{\hbox{$#1\skew\skewfactor\theaccent{}$}\vss}\hidewidth}}

\newenvironment{PROOF}[2][\proofname.]
   {\begin{proof}[#1]\renewcommand{\qedsymbol}{\textsquare$_{\rm #2}$}}
   {\end{proof}}

\begin{document}

\title {Superstable Theories and Representation}
\author {Saharon Shelah}
\address{Einstein Institute of Mathematics\\
Edmond J. Safra Campus, Givat Ram\\
The Hebrew University of Jerusalem\\
Jerusalem, 91904, Israel\\
 and \\
 Department of Mathematics\\
 Hill Center - Busch Campus \\ 
 Rutgers, The State University of New Jersey \\
 110 Frelinghuysen Road \\
 Piscataway, NJ 08854-8019 USA}
\email{shelah@math.huji.ac.il}
\urladdr{http://shelah.logic.at}
\thanks{The author thanks Alice Leonhardt for the beautiful typing.
 First typed June 3, 2013. Paper 1043}


\subjclass[2010]{Primary:03C45; Secondary:03C55}

\keywords {model theory, classification theory, stability,
  representation, superstable}


\date{\today}

\begin{abstract}
In this paper we give characterizations of the first order complete superstable 
theories, in terms of an external property called
representation. In the sense of the representation property, the mentioned class of first-order theories can be regarded as ``not very complicated". This was done for "stable" and for "$ {\aleph_0}$-stable."
Here we give a complete answer for "superstable".

\end{abstract}

\maketitle
\numberwithin{equation}{section}
\setcounter{section}{-1}
\newpage

\section {Introduction}
Our motivation to investigate the properties under consideration in
this paper comes from the following
\begin{description}
\item [{Thesis}] It is very interesting to find dividing lines and it is
a fruitful approach in investigating quite general classes of models.
A {}``natural'' dividing property {}``should'' have equivalent
internal, syntactical, and external properties. ( see \cite{Sh:E53} for more)
\end{description}

Of course, we expect the natural dividing lines will have many equivalent definitions
by internal and external properties.

The class of stable (complete first order theories) $T$ is well
known (see \cite{Sh:c}), it has many equivalent definitions by ``internal,
syntactical'' properties, such as the order property. As for external properties, one may say ``for every $ \lambda \ge  |T| $ for some model $ M $ of $ T $ we
have $ \bS(M) $ has cardinality $ > \lambda $'' is such a property (characterizing instability). 
Anyhow, the property ``not having many $ \kappa $-resplendent
models (or equivalently, having at most one in each cardinality)'' is certainly such an external property (see \cite{Sh:363}).

Here we deal with another external property, {\it representability}. This notion was a try to formalize the intuition that "the class of models of a stable first order theory is not much more complicated than the class
of models $ M=(A, \dots, E_t, \dots)_{s \in I } $
where $E^M_t$ is an equivalence relation on $A$ refining $E^M_s$ for $s < t$ ; and  $I$ is a linear order of cardinality $\le |T|$ .
 It was first defined in Cohen-Shelah \cite{CoSh:919}, where it was shown that one may characterize stability and $\aleph_0$-stability by means of representability.  In this paper we give a complete answer also for the superstable case. Moreover if $T$ is uncountable we consider other 
values of $\kappa ( T )$. That is, recall that for a stable 
(complete first order) theory $T$, $\kappa (T)$ 
can be any cardinal in the interval $ [{\aleph_0}, |T|^+ ) $. So if $T$ is countable there are two possible values- ${\aleph_0}, {\aleph_1}$,
the second is dealt with in \cite{CoSh:919} and the first
in Theorem \ref{d2}. But if $T$ is uncountable, the result
above gives a representation in a class which depends just
on $ |T| $, so it is natural to suspect that if 
$ \kappa (T ) < |T|^+ $ we can restrict this class further.
  We succeed to  do this in Theorem \ref{e8}.

The results are phrased below, and the full definition appears in Definition \ref{def:representation}, but first 
consider a simplified version.
We say that a a model $ M $ is $\mfk$-representable for a class $ {\frak k}   $ when
there exists a structure $\I \in \mfk$
with the same universe as $ M $ such that for any
$n$ and two sequences of length $n$ from $M$, if they realize
the same quantifier free
type in $ \I $ then they realize the same
(first order) type in $M$.  Of course, $T$ is
$\mfk$-representable if every model of $ T $ is $\mfk$-representable.
We prove, e.g. that $ T $ is superstable iff it is $ {\frak k}_\kappa ^{\rm unary}$-representable 
for some $ \kappa $ where $ {\frak k}_\kappa ^{\rm unary} $
is the class of structures with exactly $ \kappa $ unary functions (and nothing else).

There is also a relative characterizing ``$ \kappa (T)  \le \kappa^+$ ''.

This raises various further questions

\begin{description}
\item[{Problem}]\  
\begin{enumerate}
\item Can we characterize, by representability ``$ T $ is strongly
dependent '', similarly for the various relatives (see [Sh:863])

\item For a natural number $ n $ , what is the class
of $ T $ representable by $ \mfk_ \kappa ^n $
of structures with just $ \kappa  $ $ n $-place
functions (or relations).
\end{enumerate}
\end{description}

The main results presented in this paper are:
\begin{description}
\item [{Characterization~of~superstable~theories~(Theorem~\ref{d2})}]~
\item [{Characterization~of~$\kappa(T)$~(Theorem~\ref{e8})}]~

\end{description}

In the attempt to extend the framework of representation it seemed 
natural, initially, to conjecture that if we consider representation over linear orders 
rather than over sets, we could find an analogous characterizations for dependent theories.
However, such characterizations would imply strong theorems on existence
of indiscernible sequences. In \cite{KpSh:946}, some dependent
theories were discovered for which it is provably ``quite hard to find indiscernible
subsequences'', implying that this conjecture would fail in its original formulation.

The  author thanks Yatir Halevi for doing much to improve the presentation.

\section{Structure Classes and Representations}
We recall some needed definitions and properties from \cite{CoSh:919}.
\begin{conv}\label{conv:structure-classes}
$(1)$ The vocabulary is a set of individual constants, (partial) function symbols and relation symbols (=predicates),
each with the number of places (=arity) being finite. Individual
constants may be considered as $0$-place function symbols; here function symbols are interpreted as partial functions.

$(2)$ A structure $\I=\left\langle \tau,I,\models\right\rangle $
is a triple of vocabulary, universe(domain) and the interpretation
relation for the vocabulary: let $|\I|=I$, $\|\I\|$ the cardinality of $I$ and $\tau_\I=\tau$; $\I$ is called a $\tau$-structure.

$(3)$ $\mfk$ denotes a class of structures in a given vocabulary
$\tau_\mfk$, so $\I\in\mfk\ \Rightarrow\ \I\text{ is a }\tau_\mfk\text{-structure}$.

$(4)$ $\mathcal{L}^\tau_{\rm qf}$ denotes all the quantifier-free formulas with terms from $\tau_\mfk$.
That is, finite boolean combinations of atomic formulas, where atomic formulas (for $\tau$) have the form 
$P(\sigma_0,\ldots,\sigma_{n-1})$ or $\sigma_0=\sigma_1$ for some $n$-ary predicate $P\in \tau$, $\sigma_0\ldots$ are terms, i.e. they are in the
closure of the set of variables by function (and partial function) symbols.

$(5)$ If $\I$ a $\tau$-structure, 
$\bar{a}=\left<a_i:i<\alpha\right>\in \!^\alpha|\I|$, then \[
\tp_{\rm qf}\left(\bar{a},B,\I\right)=
\left\{\varphi(\bar{x},\bar{b}):\varphi(\bar{x},\bar{y}) \in \mathcal{L}^\tau_{\rm qf}:\I\models \varphi(\bar a,\bar b),\ \bar{b}\in \!^{\lg(\bar{y})}B\right\}
\]
\end{conv}

\subsection{Defining representations}
We recall the definition of a representation.
\begin{defn}
\label{def:representation} Consider a structure $\I$.

\begin{enumerate}
\item For a structure $\J$, a function $f:|\I|\rightarrow |\J|$ is called a representation of $\I$ in $\J$ if  \[
\tp_{\rm qf}(f(\overline{a}),\emptyset,\J)=\tp_{\rm qf}(f(\overline{b}),\emptyset,\J)\;\Rightarrow\;
\tp(\overline{a},\emptyset,\I)=\tp(\overline{b},\emptyset,\I)\]
for any two sequences $\overline{a},\overline{b}\in\!^{<\omega}I$ of equal length.

\item We say that $\I$ is represented in a class of models $\mfk$ 
if there exists a $\J\in\mfk$ such that $\I$ is represented
in $\J$.

\item For two classes of structures $\mfk_{0},\mfk$ we
say that $\mfk_{0}$ is represented in $\mfk$
if every $\I\in\mfk_{0}$ is represented in $\mfk$.

\item We say that a first-order theory $T$ is represented in
$\mfk$ if the elementary class $\EC(T)$ of $T$ is represented in $\mfk$.

\end{enumerate}

\end{defn}

\begin{defn}
\label{def:k^eq}$\keq$ denotes the class of 
structures of the vocabulary $\left\{ =\right\}$, where $\rm{eq}$ stands for equality.
\end{defn}

\subsection{The free algebras $\M_{\mu,\kappa}$}

\begin{defn}\label{def:free_algebra_M}
For a given structure $\I$, we define the structure $\M_{\mu,\kappa}(\I)$ as 
the structure whose vocabulary is $\tau_\I\cup\{ F_{\alpha,\beta}:\alpha<\mu,\beta<\kappa\}$,
with a $\beta$-ary function symbol $F_{\alpha,\beta}$ for all $\alpha<\mu,\beta<\kappa$.
(the vocabulary of $\I$ includes a unary relation symbol $I$ for the structure's universe, and we will assume 
$F_{\alpha,\beta}\notin \tau_\I$).
The universe for this structure is\footnote{This defines a set and not a proper class by remark \ref{rem:M_mu_kappa_is_a_set}.
}:\[
\M_{\mu,\kappa}(\I)=\bigcup_{\gamma\in{\rm Ord}}\M_{\mu,\kappa,\gamma}(\I)\]
 Where $\M_{\zeta}=\M_{\mu,\kappa,\zeta}(\I)$ is defined as follows:\end{defn}
\begin{itemize}
\item $\M_{0}(\I):=|\I|$
\item For limit $\zeta$: $\M_{\zeta}(\I)=\bigcup_{\xi<\zeta}\M_{\xi}(\I)$
\item For $\zeta=\gamma+1$ \[
\M_{\zeta}=\M_{\gamma}\cup\left\{ F_{\alpha,\beta}(\overline{b}):\overline{b}\in\!^{\beta}\M_{\gamma},\;\alpha<\mu,\;\beta<\kappa\right\} \]

\end{itemize}
Where $F_{\alpha,\beta}(\overline{b})$ is treated as a formal object.
The symbols in $\tau_\I$ have the same interpretation as in $\I$. In particular,
$\alpha$-ary functions may be interpreted as $(\alpha+1)-ary$ relations.
The $\beta$-ary function $F_{\alpha,\beta}(\overline{x})$ is interpreted
as the mapping $\overline{a}\mapsto F_{\alpha,\beta}(\overline{a})$
for all $\overline{a}\in\!^{\beta}\left|\M_{\mu,\kappa}(I)\right|$,
where $F_{\alpha,\beta}(\overline{a})$ on the right side of the mapping
is the formal object. If $\mu=\kappa=\aleph_0$ we may omit them.

\begin{remark}
\label{rem:M_mu_kappa_is_a_set}It is shown in \cite{CoSh:919} that
$\M_{\mu,\kappa}(S)$ is a set (though defined as a class).
\end{remark}

\subsection{Extensions of classes of structures}

\begin{disc}
For a class of structures $\mfk$, we define several classes
of structures that are based on $\mfk$.
\end{disc}

\begin{defn}
\label{def:struct_class_enrich_ex0}$\Ex_{\mu,\kappa}^{0}(\mfk)$
is the class of structures $\I^{+}$ which, for some $\I\in\mfk$ satisfy
$|\I^+|=|\I|$;$\tau_{\I^+}=\tau_{\I} \cup \left\{ P_{\alpha}:\alpha<\mu\right\} \cup 
\left\{ F_{\beta}:\beta<\kappa\right\}$ for new unary relation symbols $P_\alpha$ and new unary
function symbols $F_{\beta}$; such that if $\mu>0$ then
$\left\langle P_{\alpha}^{\I^{+}}:\alpha<\mu\right\rangle $
is a partition of $\left|\I\right|$;
and, $\left\langle F_{\beta}^{\I^{+}}:\beta<\kappa\right\rangle $
are {\bf partial} unary functions.
\end{defn}

\begin{defn}
\label{def:struct_class_enrich_ex0_lf}
$\Ex_{\mu,\kappa}^{0,{\rm lf}}(\mfk)$
is the class of structures in $\Ex_{\mu,\kappa}^{0}(\mfk)$
for which the closure of every element under the new functions is
finite. ( ${\rm lf}$ stands for {}``locally finite'').
\end{defn}

\begin{defn}\label{def:struct_class_enrich_ex1}\ 
 $\Ex_{\mu,\kappa}^{1}(\mfk)$
is the class of structures in $\Ex_{\mu,\kappa}^{0}(\mfk)$
for which $F_{\beta}(P_{\alpha})\subseteq P_{<\alpha}:=\bigcup_{\gamma<\alpha}P_{\gamma}$
holds for every $\alpha<\mu,\beta<\kappa$.

\end{defn}

\begin{defn}
\label{def:struct_class_enrich_ex2}$\Ex_{\mu,\kappa}^{2}\left(\mfk\right)$
is the class of structures of the form $\I^+=\mathcal{M}_{\mu,\kappa}(\I)$, for some $\I\in\mfk$ ( cf. Definition \ref{def:free_algebra_M} ).
\end{defn}

\begin{conv}
 $\Ex_{\mu,\kappa}$ will be one of the above classes.
\end{conv}

\section{Superstable theories}

The main theorem in this section is
\begin{theorem}
\label{d2}
For a first-order, complete theory $T$ the following are equivalent:
\mn
\begin{enumerate}
\item  $T$ is superstable.
\sn
\item  $T$ is representable in $\Ex^2_{2^{|T|},\aleph_0}(\gt^{\eq})$
\sn
\item  $T$ is representable in $\Ex^1_{2^{|T|},2}(\gt^{\eq})$
\sn
\item  $T$ is representable in $\Ex^{0,\lf}_{2^{|T|},2}(\gt^{\eq})$
\sn
\item  $T$ is representable in $\Ex^2_{\mu,\aleph_0}(\gk^{\eq})$
for some cardinal $\mu$
\sn
\item  $T$ is representable in $\Ex^{0,\lf}_{\mu,\kappa}(\gk^{\eq})$
for some cardinals $\mu,\kappa$.
\end{enumerate}
\end{theorem}

\begin{PROOF}{\ref{d2}}
$2 \Rightarrow 5, 4 \Rightarrow 6$ are immediate. 

$2 \Rightarrow 3$ is direct from \cite[1.30]{CoSh:919}

$3 \Rightarrow 4$ direct from \cite[1.24]{CoSh:919} 

$5 \Rightarrow 6$ 

This follows since $\Ex^2_{\mu,\aleph_0}(\gk^{\eq})$ is qf-representable in 
$\Ex^1_{\mu,2}(\gk^{\eq})$ by
\cite[1.30]{CoSh:919} and $\Ex^1_{\mu,2}(\gk^{\eq}) \subseteq 
\Ex^{0,\lf}_{\mu,2}(\gk^{\eq})$ by \cite[1.24]{CoSh:919} with 2 here
standing for $\kappa$ there.

The rest follows from Theorem \ref{d8} below giving
$1 \Rightarrow 2$ and Theorem \ref{d5} below giving $6 \Rightarrow 1$.
\end{PROOF}

\begin{theorem}
\label{d5}
If $T$ is representable in $\Ex_{\mu,\kappa}^{0,\lf}(\gk^{\eq})$
for some cardinals $\mu,\kappa$ then $T$ is superstable.
\end{theorem}

\begin{PROOF}{\ref{d5}}
Like the proof of Propositions \cite[Th.2.4,2.5]{CoSh:919}.
\end{PROOF}

\begin{theorem}
\label{d8}
Every superstable $T$ is representable
in $\Ex_{2^{|T|},\aleph_0}^2(\gk^{\eq})$.
\end{theorem}

\begin{PROOF}{\ref{d8}}
Let $T$ be superstable.  Let $M \prec \gC_T$. 
We choose $B_n,\langle a_{s},u_{s}:s\in S_{n}\rangle$ by 
induction on $n<\omega$ such that:
\mn
\begin{enumerate}
\item[$\circledast_0$]
\item[${{}}$] $(1)\quad $ $S_{n} \cap S_{k} = \emptyset \;(k<n)$
\sn
\item[${{}}$] $(2)\quad $  $\langle a_{s}:s \in S_{n} \rangle \subseteq M$
\sn
\item[${{}}$] $(3)\quad $  $B_{n} = \{a_{s}:s \in S_{<n}\} \subseteq M$, where 
$S_{<n} := \cup\{S_{k}:k < n\}$, as usual
\sn
\item[${{}}$] $(4)\quad $  $\langle a_{s}:s \in S_{n}\rangle$ is without repetitions, 
disjoint from $\{a_{s}:s \in S_{<n}\}$ and 

\hskip15pt independent over $B_{n}$, 
\sn
\item[${{}}$] $(5)\quad $  for all $s \in S,u_{s} \subseteq S_{<n}$ is finite such 
that $t \in u_{s} \Rightarrow u_{t} \subseteq u_{s}$
and 

\hskip15pt $\tp(a_{s},B_{n})$ does not fork over $\{a_{t}:t \in u_{s}\}$
\sn
\item[${{}}$] $(6)\quad $  $\langle a_{s}:s \in S_{n}\rangle$ is maximal under conditions
  1-5.
\end{enumerate}
\mn
Here we make a convention that $u,v,w$ vary on $\cI$ defined below:
\mn
\begin{enumerate}
\item[$\circledast_1$]   Since $T$ is superstable, it is possible to carry the induction.
\sn
\item[$\circledast_2$]  $(0) \quad \cI = \{u:u \subseteq S, u \text{ finite}\}, \text{ where } S =
\bigcup\limits_{n} S_n$
\sn
\item[${{}}$]  $(1) \quad$ for $v \in \cI$ let $\cl(v)$ be the minimal 
$u \supseteq v$ such that $u_t \subseteq u$ holds 

\hskip15pt for all $t \in u$;
\sn
\item[${{}}$]  $(2) \quad$ we define $\cI^{\cl} = \{u \in \cI:u =
  \cl(u)\}$;

\sn
\item[$\circledast_3$]  $(0) \quad$ if $u \in \cI$ then $\cl(u) \in
  \cI$
\sn
\item[${{}}$]  $(1) \quad v \subseteq u \Rightarrow \cl(v) \subseteq
  \cl(u)$;
\sn
\item[${{}}$]  $(2) \quad \cl(u_1 \cup u_2) = \cl(u_1) \cup \cl(u_2)$;
\sn
\item[${{}}$]  $(3) \quad \cl(\{s\}) = u_s \cup\{s\}$;
\sn
\item[${{}}$]  $(4) \quad \cl(\cl(u)) = \cl(u)$;
\sn
\item[${{}}$]  $(5) \quad \cl(u) = \bigcup\{u_s:s\in u\}\cup u$
\sn
\item[$\circledast_4$]   $|M| = \{a_s:s \in S\}$
(Why? Otherwise, there exists $a \in |M| 
\backslash \{a_s:s\in S\}$ such that (since $T$ is superstable) 
$\tp(a,\{a_s:s\in S\})$ does not fork over $\{a_s:s \in v\}$ for some
finite subset $v \subseteq S$.
Let $u = \cl(v)$, so $u \in \cI^\cl$ and let $n$ be such that $u
\subseteq S_n$ and we get a contradiction to the 
maximality of $\{a_s:s\in S_n\}$.)
\sn
\item[$\circledast_5$]   Let $\langle v_{\alpha}:\alpha <
  \alpha(*)\rangle$ enumerate $\cI$ (without repetition) such that
\sn
\item[${{}}$]  $(1) \quad v_{\alpha} \subseteq v_{\beta} \Rightarrow
\alpha \le \beta$;
\sn
\item[${{}}$]  $(2) \quad \alpha < \beta \wedge v_{\beta} \subseteq 
S_{<n} \Rightarrow v_{\alpha} \subseteq S_{<n}$.  

\end{enumerate}
\mn
We choose a model $M_{v_\alpha}$ by induction on $\alpha$ such that:
\mn
\begin{enumerate}
\item[$\circledast_6$]  $(1) \quad M_{v_\alpha} \prec \gC_T$ has 
cardinality $\le \aleph_{0} + |T|$;
\sn
\item[${{}}$]  $(2) \quad v_{\alpha} \subseteq v_{\beta} \Rightarrow 
M_{v_\alpha} \prec M_{v_\beta}$; 
\sn
\item[${{}}$]  $(3) \quad \bigcup\{M_{v_\beta}:\beta < \alpha \wedge 
v_\beta \subseteq v_\alpha\} \subsetneq  M_{v_\alpha}$;
\sn
\item[${{}}$]  $(4) \quad$ if $s \in v_\alpha$ and $u_s 
\subseteq v_\alpha$ then $a_s \in M_{v_\alpha}$;
\sn
\item[${{}}$]  $(5) \quad \tp(M_{v_\alpha},\cup\{M_{v_\beta}:\beta < 
\alpha\}\cup M)$ does not fork over $B_{v_\alpha} := \cup\{M_{v_\beta}:
v_{\beta}\subseteq v_{\alpha},\:\beta<\alpha\} \cup\{a_{s}:\cl(\{s\}) 
\subseteq v_{\alpha}\}$
\sn
\item[$\circledast_7$]   By $\circledast_6(3)$, clearly $\alpha <
  \beta \Rightarrow M_{v_\alpha} \ne M_{v_\beta}$.
\end{enumerate}
\mn
The induction is clearly possible.

A major point is
\mn
\begin{enumerate}
\item[$\circledast_8$]  $\tp(M_{v_\alpha},\cup\{M_{v_\beta}:\beta<\alpha\})$ does not fork over 
$A_{v_\alpha}:=\cup\{M_{v_\beta}:v_\beta
\subsetneq v_\alpha\}$.
\end{enumerate}
\mn
[Why?  If $v_\alpha=\emptyset$ this is trivial so assume 
$v_\alpha \ne \emptyset$. 

Let $n$ be such that $v_\alpha \subseteq
  S_{\le n},v_\alpha \nsubseteq S_{<n}$ and 
\begin{enumerate}
\item[$\circledast_{8.1}$] let $\langle t_\ell:\ell <
  k\rangle=\langle t^\alpha_\ell:\ell <
  k_\alpha\rangle$ list $\{s\in v_\alpha: s\notin S_{<n} \text{ and } 
\cl(\{s\})\subseteq v_\alpha\}.$
\end{enumerate}

First, assume $k=0$. So if $s \in v_{\alpha}$ and $\cl(\{s\})
\subseteq v_\alpha$ then $s \in v_\alpha\cap S_{<n}$, this implies
that $u_s\cup\{s\}=\cl(\{s\}) \subseteq S_{<n}$, hence by $\circledast_6(4)$, $a_s \in 
M_{v_\alpha \cap S_{<n}} \subseteq A_{v_\alpha}$. This implies that 
$B_{v_\alpha}\subseteq A_{v_\alpha}$ (in fact equal - see their 
definitions in $\circledast_6(5),\circledast_8$). 
Now $\circledast_6(5)$ says that
$\tp(M_{v_\alpha},\cup\{M_{v_\beta}:\beta < \alpha\})$ does not 
fork over $B_{v_\alpha}$, so by monotonicity of non-forking and 
the last sentence, it does not fork over $A_{v_\alpha}$ as desired. 

Second, assume $k=1$ and $(\forall \beta < \alpha) \left( \cl(\{t_0\})
  \nsubseteq v_\beta\right)$,

Since necessarily $v_\alpha=cl(\{t_0\})$ so $B_{v_\alpha} = A_{v_\alpha}\cup\{a_{t_0}\}$ and clearly $\tp(a_{t_0},\cup\{M_{v_\beta}:\beta<\alpha\}$ does not fork over $A_{v_\alpha}$, together with $\circledast_6(5)$ we get that $\tp(M_{v_\alpha},\cup\{M_{v_\beta}:\beta<\alpha\})$ 
does not fork over $A_{v_\alpha}$, as desired in $\circledast_8$.

Third, assume $k=1,\beta < \alpha$ and $\cl(\{t_0\}) \subseteq
v_\beta$. 
\Wilog \, $\beta$ is minimal with these properties, so necessarily 
$v_\beta=\cl(\{t_0\})$ and so again, $B_{v_\alpha} = A_{v_\alpha}$ and
we continue as in ``First'' above.

Fourth, assume $k \ge 2$. In this case, for each $\ell < k,
\cl(\{t_\ell\})$ is $v_{\beta(\ell)}$ for some unique $\beta(\ell)<\alpha$, so 
$a_{t_\ell}\in M_{v_{\beta(\ell)}} \subseteq A_{v_\alpha}$, hence, 
$B_{v_\alpha} \subseteq A_{v_\alpha}$ (in fact equal) and again
$\circledast_6(5)$ gives the desired conclusion.]

Now, $\circledast_8$ is the necessary condition in \cite[XII.2, Lemma 2.3(1)]{Sh:c} and so
we can conclude that $\langle M_v:v \in \cI\rangle$ is a stable system (see \cite[XII.2, Definition 2.1, page 598]{Sh:c}.

Now $\bar M=\langle M_{v}:v \in \cI\rangle$ is a stable system of models.
For all $v \in \cI$ let $\bar b_v$ enumerate $M_v$. By $\circledast_6(3)$, 
$\langle \bar b_v:v \in \cI \rangle$ is without repetitions.

For all $\alpha<\omega\times\omega$ we define $\cI_{\alpha}$ as follows:
\mn
\begin{enumerate}
\item[$\oplus_1$]  $(1) \quad \cI_{0} = \{\emptyset\}$;
\sn
\item[${{}}$]  $(2) \quad \cI_k = \{v \in \cI:v \subseteq
  S_{<1},|v|=k\},(k<\omega)$;
\sn
\item[${{}}$]  $(3) \quad \cI_{\omega n+k} = \{v \in \cI:v \nsubseteq
 S_{<n}, v \subseteq S_{<n+1},|v| = k+1\},(k<\omega,\:0<n<\omega)$.
\end{enumerate}
\mn
Now clearly
\mn
\begin{enumerate}
\item[$\oplus_2$]  $(1) \quad \langle \cI_\alpha : \alpha<\omega\times\omega\rangle$ is a partition of $\cI$
\sn
\item[${{}}$] $(2)\quad w \subsetneq v \in \cI_{\alpha} \Rightarrow w 
\in \cI_{<\alpha}$ for all $w,v \in \cI$
\sn
\item[$\oplus_3$] and most important,
\sn
  for all $\alpha<\omega\times\omega$ and 
$v \in \cI_{\alpha}$ the type
\[
p_v := \tp(\bar b_v,\cup\{\bar b_u:u \in \cI_{<\alpha+1} \wedge u \ne v\})
\]
does not fork over $\cup\{\bar b_w:w \subsetneq v\}$,
and moreover $p_{v}$ is the unique non-forking extension of $p_{v}
\upharpoonright \cup\{\bar b_w:w \subseteq v\}$ in 
$\mathbf{S}^{\ell g(\bar b_v)}(\cup\{\bar b_u:u \in \cI_{<\alpha+1} 
\wedge u \ne v\})$ (The proof is carried by basic properties of stable systems, see Conclusion 2.12 in \cite[Ch.XII, page 603]{Sh:c}).
\end{enumerate}

Now we define an equivalence relation $E$ on $\cI$, 
(and $E_{\alpha} = E \upharpoonright \cI_{\alpha}$) such that $v_1 
E v_2$ if and only if for some $g=g_{v_1,v_2},f=f_{v_1,v_2}$ (really $g$ 
determines $f$, we may require $g$ to be order preserving)
\mn
\begin{enumerate}
\item[$\oplus_4$]  $(1)\quad v_1,v_2 \in \cI$;
\sn
\item[${{}}$]  $(2) \quad |v_1|=|v_2|$;
\sn
\item[${{}}$]  $(3) \quad g$ is one-to-one from $v_1$ onto $v_2$ 
(we may add mapping $v_1 \cap S_{<n}$ onto $v_2 \cap S_{<n}$ for every $n$
such that if $u_i\subseteq  v_i$ (for $i=1,2$) and $g''(u_1) = u_2\Rightarrow [u_1 \in \cI_\beta 
\equiv u_2 \in \cI_\beta]$;
\sn
\item[${{}}$]  $(4) \quad u_{g(t)} = \{g(s):s\in u_t\}$, when $u_t\cup\{t\}\subseteq v_1$
\sn
\item[${{}}$]  $(5) \quad f$ is an elementary mapping of $\gC_T$
\sn
\item[${{}}$]  $(6) \quad$ if $u_\ell \subseteq v_\ell$ for 
$\ell=1,2$ and $g$ maps $u_1$ onto $u_2$, then $f$ maps 
$\bar{b}_{u_1}$ to $\bar{b}_{u_2}$
\sn
\item[${{}}$]  $(7) \quad \Dom(f) = \cup\{\bar b_u:u \subseteq v_1\}$
\item[${{}}$] $(8) \quad$ if  $s_l\in v_l,$  $cl({s_l})\subseteq v_l$ for $l=1,2$ and $g(s_1)=s_2$ then $f(a_{s_1})=a_{s_2}$.
\end{enumerate}
\mn
(element-by-element, and this implies $f_{v_{1},v_{2}}$ is unique). 
(So for some bijection $g_{v_{1},v_{2}}:v_{1}\to v_{2}$ which
preserves being in $\cI_{\beta}$ for all $\beta<\alpha$, such that
$f_{v_{1},v_{2}}$ maps $\bar b_{w_1}$ to $\bar b_{g_{v_{1},v_{2}}(w_{1})}$
for all $w_{1}\subset v_{1}$.)  Let $\langle I_{\alpha,i}:
i<i(\alpha) \le 2^{|T|}\rangle$ enumerate the equivalence classes 
of $E_{\alpha}$. 

We get that $\boxplus_1 \Rightarrow \boxplus_2$ where
\mn
\begin{enumerate}
\item[$\boxplus_1$]  $(\alpha) \quad$ the sets $\{v_{0} \ldots v_{n-1}\},
\{u_{0} \ldots u_{n-1}\} \subset \cI$ are closed under subsets
\sn
\item[${{}}$]  $(\beta) \quad \bigwedge_{\alpha,i}[v_{l} \in 
I_{\alpha,i} \Leftrightarrow u_{l} \in I_{\alpha,i}]$
\sn
\item[${{}}$]  $(\gamma) \quad u_{l(1)} \subset u_{l(2)} 
\Leftrightarrow v_{l(1)} \subset u_{l(2)}$ for all $l(1),l(2)<n$
\sn
\item[${{}}$]  $(\delta) \quad$ if $v_{l(1)} \subseteq v_{l(2)}$ 
then $g_{v_{l(2)},u_{l(2)}}$ maps $v_{l(1)}$ onto $u_{l(1)}$.
\sn
\item[$\boxplus_2$]   The sequences $\bar b_{v_0} \frown \ldots \frown
\bar b_{v_{n-1}}$ and $\bar b_{u_0} \frown \ldots \frown \bar b_{u_{n-1}}$
realize the same complete type over $\emptyset$. (This follows from
the definitions of the equivalence relations $E_{\alpha}$ and $\boxtimes_1$
above).
\end{enumerate}
\mn
%
%

Let $\mathcal{A}$ be the free algebra, i.e. a structure of the form $\mathcal{M}_{2^{|T|},\aleph_0}(\keq)$, generated by $\{x_u: u\in \cI\}$ and the function symbols $\{F_{e,\zeta}: \zeta<|T|,e=u/E \text{ for some $u\in\cI$}\}$, each one of them $|\{v\in \cI: v\subseteq u\}|$-ary for $u\in \cI$ such that $e=u/E$. Since for every $b\in M$ there is a unique $(u,\zeta)$ such that $b\in \bar b_{u,\zeta}$ (where $\bar b_{u,\zeta}$ is the $\zeta$-element of $\bar b_u$), recalling $\oplus_4(8)$, the following gives the desired representation
\[H(b)=\langle \dots x_v \dots\rangle_{v\subseteq u} \smallfrown \langle F_{u/E,\zeta} (\dots, x_v,\dots)_{v\subseteq u}\rangle,\]
after choosing some ordering on the elements of $\cI$.
\end{PROOF}

\section {Between stable and superstable}

\begin{discussion}
\label{e3}
For superstable $T$, we may wonder about whether ``the cardinal $2^{|T|}$ is optimal". 
Really, $\lambda(T)$ is sufficient where
\mn
\begin{enumerate}
\item[$(*)_{1.1}$]  $\lambda(T) = \min\{\lambda:T \text{ is stable in } 
\lambda\}$.
\end{enumerate}
\mn                                     
Note that 
\mn
\begin{enumerate}
\item[$(*)_{1.2}$]  If $T$ is countable then $\lambda(T) = \aleph_0$
is equivalent to $T$ is $\aleph_0$-stable and 
\sn
\item[$(*)_{1.3}$]  if $T$ is countable and $\lambda(T) > \aleph_0$
  then $\lambda(T)=2^{\aleph_0}$.
\end{enumerate}
\mn
\end{discussion}

\begin{theorem}
\label{e6}
In Theorem \ref{d8}, $\Ex^2_{\lambda(T),\aleph_0}(\gk^{\eq})$ suffice.
\end{theorem}
\begin{PROOF}{\ref{e6}}
We repeat the proof of Theorem \ref{d8} with some changes. Choosing $M_{v_\alpha}$ by induction on $\alpha$ we now demand more than in $\circledast_6$ in the proof of Theorem \ref{d8}:
\mn
\begin{enumerate}
\item[$\circledast_6^*$]  $(1) \quad M_{v_\alpha} \prec \gC$ is saturated 
of cardinality $\lambda(T)$;
\sn
\item[${{}}$]  $(2) \quad v_\beta \subseteq v_\alpha 
\Rightarrow M_{v_\beta} \subseteq M_{v_{\alpha}}$
\sn
\item[${{}}$]  $(3) \quad (M_\alpha,c)_{c \in \cup \{M_{v_\beta}:
v_\beta \subset v_\alpha\}}$ is saturated.
\end{enumerate}
\mn

In $\circledast_{8.1}$ of the proof of Theorem \ref{d8} we replace $v_\alpha$ with $v$, i.e. 
\begin{enumerate}
\item[$\circledast_{8.1}^*$] if $v\in S_n$ then $\bar t_v=\langle t_{v,\ell}:\ell<k_v\rangle$ list $\{s\in v:s\notin S_{<n} \text{ and } \cl (\{s\})\subseteq v\}$ in increasing order.
\end{enumerate}

Now after $\oplus_3$ in the proof and as in the proof, we choose by induction on $\alpha<\omega\times\omega$,  $E_\alpha, f_{v_1,v_2},g_{v_1,v_2}$ for $(v_1,v_2)\in E_\alpha$ such that
\begin{enumerate}
\item[$\oplus_4^*$] 
\begin{enumerate}
\item[(a)] $E_\alpha$ is an equivalence relation on $\cI_\alpha$ with $\leq \lambda (T)$ equivalence classes
\item[(b)] Clauses $(1)-(8)$ of $\oplus_4$ holds with the following modifications: if $(v_1,v_2)\in E_\alpha$ then $f=f_{v_1,v_2},\, g=g_{v_1,v_2}$ satisfy
\sn  

\item[${{}}$] $(1)\quad v_1,v_2 \in \cI$;
\sn
\item[${{}}$]  $(2) \quad |v_1|=|v_2|$;
\sn
\item[${{}}$]  $(3) \quad g$ is one-to-one from $v_1$ onto $v_2$ 
(we may add mapping $v_1 \cap S_{<n}$ onto $v_2 \cap S_{<n}$ for every $n$
such that if $u_i\subseteq  v_i$ (for $i=1,2$) and $g''(u_1) = u_2\Rightarrow [u_1 \in \cI_\beta 
\equiv u_2 \in \cI_\beta]\wedge g\restriction u_1=g_{u_1,u_2}$;
\sn
\item[${{}}$]  $(4) \quad u_{g(t)} = \{g(s):s\in u_t\}$, when $u_t\cup\{t\}\subseteq v_1$
\sn
\item[${{}}$]  $(5) \quad f_{v_1,v_2}$ is an elementary mapping of $\gC_T$
\sn
\item[${{}}$]  $(6) \quad$ if $u_\ell \subseteq v_\ell$ for 
$\ell=1,2$ and $g$ maps $u_1$ onto $u_2$, then $f_{v_1,v_2}$ maps 
$M_{u_1}$ onto $M_{u_2}$ and $f\restriction M_{u_1}=f_{u_1,u_2}$
\sn
\item[${{}}$]  $(7) \quad \Dom(f_{v_1,v_2}) = \cup\{M_u:u \subseteq v_1\}$
\item[${{}}$] $(8) \quad$ if  $s_l\in v_l,$  $cl({s_l})\subseteq v_l$ for $l=1,2$ and $g_{v_1,v_2}(s_1)=s_2$ then $f_{v_1,v_2}(a_{s_1})=a_{s_2}$.

\end{enumerate}
\end{enumerate}
Note that the $\bar b_v$ notation does not appear.

The induction is carried in the following way. For $\alpha=0$ this is obvious, so let $\alpha\in (0,\omega\times\omega)$. First we define a two place relation $E_\alpha '$ on $\cI_\alpha$:
\begin{enumerate}
\item[$\oplus_{5,\alpha}$] $v_1E_\alpha' v_2$ if and only if
\begin{enumerate}
\item[$(1)$] $v_1,v_2\in \cI_\alpha$ and $|v_1|=|v_2|$
\item[$(2)$] letting $v_\ell=\{t_{\ell,j}:j<|v_\ell|\}$ for $\ell=1,2$ list $v_\ell$ in increasing order we have
\[j<|v_1|\Rightarrow \bigvee_{\beta <\alpha} (v_1\setminus \{t_{1,j}\}) E_\beta (v_2\setminus \{t_{2,j}\})\]
\item[$(3)$] $g=g_{v_1,v_2}$ maps $t_{1,j}$ to $t_{2,j}$ for $j<|v_1|$.
\end{enumerate}
\end{enumerate}

Note that for every $\alpha$, $E_\alpha'$ is an equivalence relation on $\cI_\alpha$ with $\leq \lambda (T)$ equivalence classes. By \cite[XII, 3.2 page 604, 3.5 page 608]{Sh:c} and the induction hypothesis, if $(v_1,v_2)\in E_\alpha'$ then $f_{v_1,v_2}':=\bigcup\{f_{v_1\setminus \{t_{1,j}\},v_2\setminus\{t_{2,j}\}}:j<|v_1|\}$ is an elementary map.

\begin{enumerate}
\item[$\oplus_{7,\alpha}$] We define the place relation $E_\alpha$ on $\cI_\alpha$, as $v_1E_\alpha v_2$ if and only if $v_1E_\alpha' v_2$ and $f_{v_1,v_2}$ and $f_{v_1,v_2}'$ maps $\tp(\langle a_{t_{v_1,j}}:j<|v_1|\rangle,\bigcup \{M_u :u\subsetneq v_1\})$ onto $\tp(\langle a_{t_{v_2,j}}:j<|v_2|\rangle,\bigcup \{M_u :u\subsetneq v_2\})$.
\end{enumerate}

Since $T$ is stable in $\lambda(T)$, $E_\alpha$ is an equivalence relation on $\cI_\alpha$ with $\leq \lambda(T)$ equivalence classes refining $E_\alpha'$.

Finally, if $(v_1,v_2)\in E_\alpha$ then let $f_{v_1,v_2}$ be an elementary map from $M_{v_1}$ to $M_{v_2}$ extending $f_{v_1,v_2}'$ and mapping $\langle a_{t_{v_1,j}}:j<|v_1|\rangle$ to $\langle a_{t_{v_2,j}}:j<|v_2|\rangle$

Now that $E_\alpha, f_{v_1,v_2}$ and $g_{v_1,v_2}$ are defined inductively, the rest of the proof is as in the proof of Theorem \ref{d8}.
\end{PROOF}

\bigskip

\subsection {Characterization of $\kappa(T)$} \
\bigskip

Let $T$ be a stable theory. Recall the definition of $\kappa (T)$ from \cite[III.3]{Sh:c}:
$\kappa (T)=\sup_m \kappa^m(T)$, where $\kappa^m(T)$ is the first infinite $\kappa$ for which we do not have an increasing sequence $A_i$, $i\leq \kappa$ and $p\in S_m(A_\kappa)$ such that for all $i<\kappa$, $p \upharpoonright A_{i+1}$ forks over $A_i$. 

\begin{remark}
If $T$ is stable then $\kappa (T)\leq |T|^+$ \cite[Corollary 3.3]{Sh:c}. $T$ is superstable if and only if  $\kappa(T)=\aleph_0$ \cite[Corollary 3.8]{Sh:c}.
\end{remark}

\begin{definition}\label{f5}
Assume $T$ is stable. \begin{enumerate}
\item We say that $\bold a$ is an independent system 
where $\bold a=\langle \cI_{\bold a},\bar A_{\bold a}\rangle$ consists of 
\sn
\begin{enumerate}
 
\item[$\bullet_1$]   $\bar A_{\bold a} = \langle A_{\bold a,v}:v \in \cI_{\bold a} \rangle$;
\sn
\item[$\bullet_2$]   $\cI_{\bold a}$ is a family of subsets of $S=\cup \cI_{\bold a}$;
\sn
\item[$\bullet_3$]  $\emptyset \in \cI_{\bold a}$;
\sn
\item[$\bullet_4$]  $\cI_{\bold a}$ is closed under any intersection
\sn
\item[$\bullet_5$]  $v \subseteq u \Rightarrow A_{\bold a,v} \subseteq A_{\bold a,u}$ 
for all $v,u \in \cI_{\bold a}$
\sn
\item[$\bullet_6$]  $A_{\bold a,v \cap u} = A_{\bold a,v} \cap A_{\bold a,u}$ for all $v,u \in \cI_{\bold a}$
\sn
\item[$\bullet_7$]  $\tp(A_{\bold a,v},\cup\{A_{\bold a ,u_\ell}:\ell<n\}$ does not fork 
over $\cup\{A_{\bold a,v\cap u_\ell}:\ell<n\}$ for all $v,u_0,\ldots,u_{n-1}\in \cI_{\bold a}$

\end{enumerate}
\item We may also say that $\bold a$ is an independent $\cI$-system and let $\cI_{\bold a}=\cI$ etc. \\
We say that $\bold a$ is an independent $(\theta,\cI)$-system if in addition we demand that $|A_v|\leq \theta$ for all $v\in \cI$.

\item Let $\bold a \bold b$ mean:
\begin{enumerate}
\item[$\bullet$] $\bold a, \bold b$ are independent systems,
\item[$\bullet$] $\cI_{\bold a}=\cI_{\bold b}$,
\item[$\bullet$] $A_{\bold a,v}\subseteq A_{\bold b,v}$ for $v\in \cI_{\bold a}$,
\item[$\bullet$] $\tp(A_{\bold b,v},\cup\{A_{\bold a,u}:u\in \cI_{\bold a}\})$ does not fork over $A_{\bold a,v}$.

\end{enumerate}
\end{enumerate}
\end{definition}

\begin{proposition}\label{3.5}
Assume $T$ is stable.
\begin{enumerate}
\item $\leq$ is a partial order on the the class of independent systems.
\item If $\langle \bold a_{\alpha}:\alpha<\delta\rangle$ is a $\leq$-increasing sequence then it has a least upper bound $\bold a_\delta:=\bigcup_{\alpha<\delta} \bold a_\alpha$, i.e.
\begin{enumerate}
\item[$\bullet$] $\cI_{\bold a_\delta}=\cI_{\bold a_\alpha}$ for any $\alpha<\delta$,
\item[$\bullet$] $A_{\bold a_\delta,v}=\bigcup\{A_{\bold a_\alpha,v}:\alpha<\delta\}$ for any $v\in \cI_{\bold a_\delta}$.
\end{enumerate}
\item If $\langle \bold a_\alpha:\alpha<\delta\rangle$ is a $\leq$-increasing sequence of $(\theta,\cI)$-systems and $\delta<\theta^+$ then $\bigcup	_{\alpha<\delta}\bold a_\alpha$ is an independent $(\theta,\cI)$-system.
\end{enumerate}
\end{proposition}
\begin{PROOF}{\ref{3.5}}
Should be clear.
\end{PROOF}

\begin{proposition}\label{3.6} Assume $T$ is stable and $\theta\geq |T|$. 
Let $\bold a$ be an $(\theta,\cI)$-system, with $2^{|v|}\leq \theta$ for all $v\in\cI$.
\begin{enumerate}
\item There is $(\theta,\cI)$-system $\bold b$ with $\bold a\leq \bold b$ and $\bar M=\langle M_v:v\in I\rangle$ such that for every $v\in\cI$ we have

\item[$\bullet_1$]  $A_{\bold a,v}\subseteq M_v\subseteq A_{\bold b,v}$
\item[$\bullet_2$] $M_v\prec \gC_t$
\item[$\bullet_3$] $(M_v,c)_{c\in\bigcup\{A_{\bold a,u}:u\subsetneq v,\, u\in \cI\}}$ is saturated if $T$ is stable in $\theta$.

\end{enumerate}

\end{proposition}
\begin{PROOF}{\ref{3.6}}

 Let $\langle v_\alpha :\alpha<\alpha (*)\rangle$ list $\cI$ with no repetition. We chose $M_\alpha$ by induction on $\alpha<\alpha (*)$ such that
\begin{enumerate}
\item[$(*)_1^\alpha$]\sn
\begin{enumerate}
\item[(a)] $M_\alpha\prec \gC_T$ has cardinality at most  $\theta$
\item[(b)] $A_{\bold a ,v_\alpha}\subseteq M_\alpha$
\item[(c)] $\tp(M_\alpha,\bigcup\{M_\beta:\beta<\alpha\}\cup\bigcup\{A_{\bold a,u}:u\in \cI\})$ does not for over $A_{\bold a,v}$, and 
\item[(d)] $(M_\alpha,c)_{c\in\bigcup\{A_{v_\beta}:\beta<\alpha,\, v_\beta\subsetneq v_\alpha\}}$ is saturated if $T$ is stable in $\theta$.
\end{enumerate}
\end{enumerate}

Now define $\bold b$ by setting $\bar A_{\bold b}=\langle A_{\bold b,v_\alpha}:v_\alpha\in \cI\rangle$, where $A_{\bold b,v_\alpha}:=\bigcup\{M_u:u\subseteq v_\alpha,\, u\in \cI\}$.

Note that by the non-forking calculus
\begin{enumerate}
\item[$(*)_2$]\sn
\begin{enumerate}
\item[(a)] $\bold b$ is an independent $(\theta,\cI)$-system and
\item[(b)] $\bold a\leq \bold b$.
\end{enumerate}
\end{enumerate}
%
\end{PROOF}

\begin{theorem}
\label{e8}
For a complete stable theory $T$, and $\kappa$ the following are equivalent, letting $\theta=2^{|T|}$:
\mn
\begin{enumerate}
\item   $\kappa(T) \le \kappa^+$
\sn
\item   $T$ is representable in $\Ex^2_{\theta,\kappa}(\gk^{\eq})$
\sn
\item  $T$ is representable in $\Ex^2_{\mu,\kappa}(\gk^{\eq})$
for some $\mu$
\sn
\item  $T$ is representable in $\Ex^1_{\theta,\kappa}(\gk^{\eq})$
\sn
\item  $T$ is representable in $\Ex^1_{\mu,\kappa}(\gk^{\eq})$ for some $\mu$.
\end{enumerate}
\end{theorem}

\begin{remark}
\label{e11}
In Theroem \ref{e8}, demanding $\kappa$ be regular, we can change
clauses (2)-(5) above to (2')-(5') which mean that
we use only $\bold I$ such that the closure of any set of 
cardinality $<\kappa$ has cardinality $<\kappa$ and change (1) to 
(1') $\kappa(T) \le \kappa=\delta$. The proof is similar, and in the proof we have only $S_i$ for $i<\kappa$.
\end{remark}

\begin{PROOF}{\ref{e8}}
By Theorem \ref{d2}, \wilog \, $\kappa>|T|+\aleph_0$. Let \[\partial=\min \{\lambda: \lambda \text{ regular, } \lambda\geq\kappa(T)\}.\] 

The proof is continuing the proof of the stable case and the superstable case.
So, cannibalizing the proof of Theorem \ref{d8} and/or
\cite[2.15]{CoSh:919} we have $\langle S_i:i<\partial\rangle$ is a sequence of pairwise disjoint sets of ordinals, $\langle a_s:
s \in S_i \rangle$ for $i < \partial,u_s \in [S_{<i}]^{< \kappa}$ for 
$s \in S_i,\cl(u) = \cup\{u_s:s\in u\} \cup u$ for 
$u \subseteq S=S_{<\partial}$ and $\cI = \{u:u \subseteq u_s \cup\{s\} 
\text{ for some } s \in S\}$.

By the non-forking calculus, we have that $A_{\bold a,v}=\{a_s:s\in S\text { and } \cl(\{s\})\subseteq v\}$ for $v\in\cI$ defines a $(\partial,\cI)$-independent system hence a $(\theta,\cI)$-independent system. 

Similar to the proof of Theorem \ref{d8}, we build an equivalence relation on $\cI$, but now we have to do it gradually and simultaneously with defining a new independent system. 

\begin{enumerate}
\item[$(*)_1$] Let $E_\bullet$ be the following equivalence relation on $\cI$: $v_1E_\bullet v_2$ if and only if
\begin{enumerate}
\item[(a)] The order-type of $v_1$ is equal to the order-type over $v_2$, recalling that $v_1$ and $v_2$ are set of ordinals.
\item[(b)] $g_{v_1,v_2}$, the unique order preserving map from $v_1$ onto $v_2$, maps $\{u\in \cI:u\subseteq v_1\}$ onto $\{u\in \cI: u\subseteq v_2\}.$
\end{enumerate}
\end{enumerate}

Note that $E_\bullet$ is an equivalence relation on $\cI$ with at most $\theta$ equivalence classes. Let $\bar g=\langle g_{v_1,v_2}: v_1E_\bullet v_2\rangle$.

By using Proposition \ref{3.6} repeatedly, we choose by induction on $n< \omega$  $(\bold b_n,E_n,\bar f_n,\bar M_n)$ satisfying:
\begin{enumerate}
\item[$(*)_2$]
\begin{enumerate}
\item[(a)] $\bold b_n$ is a $(\theta,\cI)$-independent system
\item[(b)] 
\begin{enumerate}
\item[(i)] $E_n$ is an equivalence relation on $\cI$ with $\leq \theta$ equivalence classes refining $E_\bullet$ and $E_m$ for $m<n$.
\item[(ii)] Let $\bar f_n=\langle f^n_{v1,v_2}:v_1E_nv_2\rangle$, where 
\item[(iii)] $f^n_{v_1,v_2}$ is an elementary mapping from $A_{\bold b_n,v_1}$ onto $A_{\bold b_n,v_2}$
\item[(iv)] if $v_2,v_3\in v_1/E_n$ then $f^n_{v_2,v_3}\circ f^n_{v_1,v_2}=f^n_{v_1,v_3}$
\item[(v)] if $cl(\{s_1\}\in v_1$ and $g^n_{v_1,v_2}=s_2$ then $f^n_{v_1,v_2}(a_{s_1})=a_{s_2}$.
\end{enumerate}

\item[(c)] if $n=m+1$ then by applying Proposition \ref{3.6}, we get a $(\theta,\cI)$-independent system $\bold b_m\leq \bold b_{n}$ $\bar M_n=\langle M_{n,v}:v\in \cI\rangle$ satisfying for every $v\in \cI$ $A_{\bold b_m,v}\subseteq M_{n,v}\subseteq A_{\bold b_n,v}$. Note that $(M_{n,v},c)_{c\in\cup\{M_{n,u}:u\subsetneq v,u\in \cI\}}$ is saturated.
\item[(d)] if $v_1E_nv_2$ for $n=m+1$ then $f^m_{v_1,v_2}\subseteq f^n_{v_1,v_2}$

\end{enumerate}
\end{enumerate}
By Proposition \ref{3.5}, $\bold b:=\bigcup_m \bold b_m$ is a $(\theta,\cI)$-independent system and $\bold a\leq \bold b$.

Clearly, for every $v\in\cI$,
\[A_{\bold b,v}=\bigcup\{A_{\bold b_m,v}:m<\omega\}\subseteq \bigcup\{M_{m,v}:m<\omega\}\]
\[\subseteq \bigcup \{A_{\bold b_{m+1}},v:m<\omega\}=A_{\bold b,v}.\] 
Hence $\gC\restriction A_{\bold b,v}\prec \gC$. 

Let $E=\cap_m E_m$ and $f_{v_1,v_2}=\bigcup_m f^m_{v_1,v_2}$. The rest is as in the proof of Theorem \ref{d8}.

\end{PROOF}

\bibliographystyle{alpha}
\bibliography{1043}

\begin{thebibliography}{She90}

\bibitem[KS]{KpSh:946}
Itay Kaplan and Saharon Shelah.
\newblock Dependent examples.
\newblock {\em preprint}.

\bibitem[SC16]{CoSh:919}
Saharon Shelah and Moran Cohen.
\newblock Stable theories and representation over sets.
\newblock {\em MLQ Math. Log. Q.}, 62(3):140--154, 2016.

\bibitem[Shea]{Sh:E53}
Saharon Shelah.
\newblock Introduction to: classification theory for abstract elementary class.
\newblock {\em TBA}.
\newblock math.LO/0903.3428.

\bibitem[Sheb]{Sh:363}
Saharon Shelah.
\newblock On spectrum of $\kappa $-respondent models.
\newblock {\em preprint}.
\newblock A revised version will appear in ``Non structure theory".

\bibitem[She90]{Sh:c}
Saharon Shelah.
\newblock {\em Classification theory and the number of nonisomorphic models},
  volume~92 of {\em Studies in Logic and the Foundations of Mathematics}.
\newblock North-Holland Publishing Co., Amsterdam, xxxiv+705 pp, 1990.

\end{thebibliography}

\end{document}